\documentclass[12pt,twoside,a4paper,reqno]{amsart}

\usepackage{amssymb}
\usepackage{eucal}
\usepackage{amsmath}

\usepackage[utf8]{inputenc}
\usepackage[T1]{fontenc}
\usepackage{amsmath,amssymb,amsthm}
\usepackage{lmodern}
\usepackage[hidelinks]{hyperref}

\newtheorem{theorem}{Theorem}[section]
\newtheorem{lemma}[theorem]{Lemma}
\newtheorem{proposition}[theorem]{Proposition}

\theoremstyle{definition}
\newtheorem{definition}[theorem]{Definition}

\theoremstyle{remark}
\newtheorem{remark}[theorem]{Remark}

\newcommand{\OO}{\mathcal O}
\newcommand{\DD}{\mathcal D}
\newcommand{\ZZ}{\mathbb Z}
\newcommand{\CC}{\mathbb C}

\newcommand{\End}{\operatorname{End}}
\newcommand{\Pic}{\operatorname{Pic}}
\newcommand{\Norm}{\operatorname{N}}

\title{Counting Complex Elliptic Curves with Self-Quotients via Class Numbers}

\author{Radu Gaba}
\address{Institute of Mathematics "Simion Stoilow" of the Romanian Academy, P.O. BOX 1-764 RO-014700 and The National University of Science and Technology Politehnica Bucharest, P.O. BOX  RO-060042
Bucharest  e-mail: {\tt radu.gaba@imar.ro}}

\begin{document}

%
%
%
%

%
%



\maketitle

\begin{abstract}

Let \(n>1\) be an integer. We consider complex elliptic curves
\(E\) admitting a finite subgroup \(C\subset E(\mathbb C)\) of
order \(n\) such that \(E/C\) is isomorphic to \(E\).
 We establish class-number counting formulas for the cyclic, noncyclic, and general cases.

More precisely, we associate to every negative quadratic-order
discriminant \(D\) the norm form
\[
Q_D(x,y)=x^2+\varepsilon_Dxy+q_Dy^2,
\]
where
\[
\varepsilon_D\in\{0,1\},\qquad
q_D=\frac{\varepsilon_D^2-D}{4}.
\]
The invariant factors of the kernel of the corresponding
endomorphism are determined by \(\gcd(x,y)\). This leads to
finite sets of admissible discriminants describing the cyclic,
noncyclic, and general cases with unequal invariant factors.

Our main result expresses the numbers of isomorphism classes in the three cases as sums of the proper ideal class numbers  \(h(D)\)  over explicitly defined finite sets of admissible discriminants.
 This provides an arithmetic class-number counting formulation of the CM classification developed in earlier work by \cite{CanepaGaba2014,CanepaGaba2019,CanepaGabaOlteanu2022}  and its subsequent algorithmic extensions.

We also characterize the equal-invariant-factor
case: it occurs precisely when $n=m^2$, and in that case
every complex elliptic curve satisfies the required condition,
with the unique subgroup being $E[m]$.

 We derive an explicit bound for the admissible discriminants, implement the counting procedure in Magma, and give numerical examples.
The proof relies on the Smith normal form of the multiplication matrix associated with a quadratic integer.

The results provide an arithmetic class-number counting
formulation of the CM classification developed in earlier work.


\end{abstract}

\section{Introduction}

Let \(E\) be a complex elliptic curve and let
\(C\subset E(\mathbb C)\) be a finite subgroup. The quotient
\(E/C\) is again a complex elliptic curve, and the natural
projection
\[
\pi:E\longrightarrow E/C
\]
is an isogeny of degree \(|C|\).

A natural problem is to determine the complex elliptic curves
\(E\) for which there exists a subgroup \(C\) of a prescribed
order \(n\) such that
\[
E/C\simeq E.
\tag{1.1}
\]

This problem has been studied through the description of
endomorphisms of complex elliptic curves and through algorithms
for the corresponding classification problem. In particular,
the CM classification and its algorithmic treatment were
developed in earlier work; see
\cite{CanepaGaba2014,CanepaGaba2019,CanepaGabaOlteanu2022}.

The purpose of the present paper is not to propose a new
classification of these elliptic curves. Instead, we give a
different arithmetic description of the counting problem.
The central observation is that the relevant elliptic curves
can be grouped according to their imaginary quadratic
endomorphism orders, and that the number of isomorphism
classes belonging to a fixed order is its proper ideal class
number.

The main result can be summarized as
\[
\boxed{
\text{number of elliptic-curve classes}
=
\text{sum of class numbers of admissible quadratic orders}.
}
\tag{1.2}
\]

The admissible orders are described explicitly by integral
solutions of a binary quadratic norm equation. Moreover,
the invariant factors of the subgroup \(C\) are determined
by the greatest common divisor of the two coordinates in
that equation.

We distinguish three situations:
\begin{enumerate}
    \item the subgroup \(C\) is cyclic;
    \item the subgroup \(C\) is noncyclic and has unequal
    invariant factors;
    \item the subgroup \(C\) has unequal invariant factors,
    without imposing cyclicity.
\end{enumerate}

The distinction between the second and third situations is
important. A single elliptic curve may admit both a cyclic
and a noncyclic subgroup of the prescribed order satisfying
(1.1). Thus the corresponding sets of discriminants may
overlap, and the general count must count each elliptic-curve
class only once.

There is also a fourth possibility, namely that the two
invariant factors of $C$ are equal. This case is fundamentally
different from the preceding ones and does not require complex
multiplication. We prove that it occurs if and only if
$n=m^2$ for some integer $m>1$. Moreover, when $n=m^2$,
every complex elliptic curve admits exactly one subgroup with
equal invariant factors satisfying (1.1), namely its
$m$-torsion subgroup $E[m]$. Consequently, the number of
isomorphism classes in this case is uncountably infinite.
This explains why the equal-invariant-factor case must be
treated separately from the finite class-number sums.

\begin{remark}
It is instructive to compare Theorem~\ref{thm:main} with the classical modular polynomial perspective. Over $\mathbb{C}$, pairs $(E, C)$ consisting of an elliptic curve $E$ together with a cyclic subgroup $C \subset E$ of order $n$ satisfying $E/C \simeq E$ correspond to roots of the diagonal modular polynomial $\Phi_n(X, X) = 0$. By Kronecker's class-number relations, the degree of $\Phi_n(X, X)$ (counting multiplicities and automorphism weights) is expressed as a sum of Hurwitz class numbers. In contrast, our formulas count \emph{isomorphism classes of elliptic curves} $E$ without multiplicity—each curve is counted at most once regardless of how many self-quotient subgroups of order $n$ it admits—and stratify the count explicitly according to the invariant factors of $C$.
\end{remark}

The paper is organized as follows. Section 2 recalls the
quadratic-order notation and the CM parametrization by proper
ideal classes. Section 3 establishes the invariant-factor
lemma using the Smith normal form. Section 4 defines the
admissible discriminants and proves the class-number counting
theorem.  Subsection 4.1 treats the equal-invariant-
factor case. Section 5 gives a finite bound for the discriminants
and illustrates the formulas with numerical examples.

\section{Quadratic orders and complex multiplication}

Let \(D<0\) be the discriminant of an imaginary quadratic
order. Thus
\[
D\equiv 0\ \text{or}\ 1\pmod 4.
\]

Put
\[
\varepsilon_D=
\begin{cases}
0,&D\equiv0\pmod4,\\
1,&D\equiv1\pmod4,
\end{cases}
\qquad
q_D=\frac{\varepsilon_D^2-D}{4}.
\tag{2.1}
\]

The corresponding imaginary quadratic order is
\[
\OO_D=\ZZ[\omega_D],
\qquad
\omega_D=\frac{\varepsilon_D+\sqrt D}{2}.
\tag{2.2}
\]
In particular,
\[
\omega_D^2
=
\varepsilon_D\omega_D-q_D.
\tag{2.3}
\]

For integers \(x,y\), the norm of
\[
\alpha=x+y\omega_D
\]
is
\[
\begin{aligned}
\Norm(\alpha)
&=(x+y\omega_D)(x+y\overline{\omega_D})\\
&=x^2+\varepsilon_Dxy+q_Dy^2.
\end{aligned}
\]

We therefore define the positive definite binary quadratic
form
\[
Q_D(x,y)
=
x^2+\varepsilon_Dxy+q_Dy^2.
\tag{2.4}
\]

Its discriminant is
\[
\varepsilon_D^2-4q_D=D.
\]

We write \(h(D)\) for the proper ideal class number of
\(\OO_D\), that is,
\[
h(D)=|\Pic(\OO_D)|.
\]
Equivalently, \(h(D)\) is the number of proper equivalence
classes of primitive positive definite binary quadratic forms
of discriminant \(D\).

We recall the standard parametrization from the theory of
complex multiplication.

\begin{proposition}[CM parametrization]
\label{prop:CM}
Let \(\OO_D\) be an imaginary quadratic order. The isomorphism
classes of complex elliptic curves \(E\) satisfying
\[
\End(E)\simeq\OO_D
\]
are in bijection with the proper ideal classes of
\(\OO_D\). Consequently, the number of such isomorphism
classes is
\[
h(D).
\]
\end{proposition}

\begin{proof}
Every complex elliptic curve can be written analytically as
\[
E=\CC/I,
\]
where \(I\subset\CC\) is a lattice. Its endomorphism ring is
\[
\End(E)=
\{\alpha\in\CC:\alpha I\subseteq I\}.
\]

The condition \(\End(E)=\OO_D\) is equivalent to the
corresponding lattice being homothetic to a proper
\(\OO_D\)-ideal. Two such lattices define isomorphic complex
elliptic curves precisely when they are homothetic.

Thus isomorphism classes of elliptic curves with
endomorphism ring \(\OO_D\) correspond to proper ideal
classes of \(\OO_D\). The number of these classes is, by
definition, \(h(D)\).

For this classical CM correspondence, see
\cite[Chapter II]{SilvermanAdvancedTopics} and
\cite[Lecture 18, Section 18.4]{SutherlandCM}.
\end{proof}

\begin{remark}
The proposition counts isomorphism classes of elliptic curves,
not pairs \((E,C)\). This distinction will be maintained
throughout the paper. In particular, if one elliptic curve
admits several subgroups satisfying (1.1), it is counted only
once in the corresponding class-number sum.
\end{remark}

\section{Invariant factors and the Smith normal form}

We now establish the arithmetic relation between an
endomorphism and the invariant factors of its kernel.

Let
\[
\alpha=x+y\omega_D\in\OO_D,
\qquad x,y\in\ZZ,
\]
and suppose that
\[
\Norm(\alpha)=n>0.
\tag{3.1}
\]

The multiplication map
\[
m_\alpha:\OO_D\longrightarrow\OO_D,
\qquad
z\longmapsto\alpha z
\]
is a \(\ZZ\)-linear map. We use the basis
\(\{1,\omega_D\}\) of \(\OO_D\).

\begin{lemma}[Invariant factors of the kernel]
\label{lem:invariant-factors}
Let \(D<0\) be a quadratic-order discriminant and let
\[
\OO_D=\ZZ[\omega_D],
\qquad
\omega_D^2=\varepsilon_D\omega_D-q_D.
\]
Suppose that
\[
\alpha=x+y\omega_D\in\OO_D,
\qquad
\Norm(\alpha)=n>0,
\]
and put
\[
g=\gcd(x,y)>0.
\]
Then
\[
\boxed{
\OO_D/\alpha\OO_D
\simeq
\ZZ/g\ZZ
\times
\ZZ/(n/g)\ZZ.
}
\tag{3.2}
\]
In particular, the invariant factors of this finite abelian
group are \(g\) and \(n/g\), and they are unequal if and only if
\[
g<\sqrt n.
\]
\end{lemma}

\begin{proof}
We use the \(\ZZ\)-basis \(\{1,\omega_D\}\) of \(\OO_D\).
Multiplication by \(\alpha=x+y\omega_D\) gives
\[
\alpha\cdot1=x+y\omega_D.
\]
Moreover, by (2.3),
\[
\begin{aligned}
\alpha\omega_D
&=x\omega_D+y\omega_D^2\\
&=x\omega_D+y(\varepsilon_D\omega_D-q_D)\\
&=-q_Dy+(x+\varepsilon_Dy)\omega_D.
\end{aligned}
\]

Consequently, the matrix of \(m_\alpha\), with respect to the
basis \(\{1,\omega_D\}\), is
\[
M_\alpha=
\begin{pmatrix}
x&-q_Dy\\
y&x+\varepsilon_Dy
\end{pmatrix}.
\tag{3.3}
\]

Its determinant is
\[
\begin{aligned}
\det M_\alpha
&=x(x+\varepsilon_Dy)+q_Dy^2\\
&=x^2+\varepsilon_Dxy+q_Dy^2\\
&=\Norm(\alpha)\\
&=n.
\end{aligned}
\tag{3.4}
\]
Thus \(M_\alpha\) is nonsingular, and its cokernel is finite.

We recall that the Smith normal form of a nonsingular
\(2\times2\) integer matrix \(M\) is a diagonal matrix
\[
\begin{pmatrix}
d_1&0\\
0&d_2
\end{pmatrix},
\qquad
d_1,d_2>0,\quad d_1\mid d_2,
\]
obtained by multiplication on the left and right by
unimodular integer matrices. The diagonal entries satisfy
\[
d_1=\gcd(\text{all entries of }M),
\qquad
d_1d_2=|\det M|.
\tag{3.5}
\]

For the matrix \(M_\alpha\), we have
\[
\begin{aligned}
\gcd(x,-q_Dy,y,x+\varepsilon_Dy)
&=\gcd(x,y)\\
&=g.
\end{aligned}
\tag{3.6}
\]
Indeed, every common divisor of \(x\) and \(y\) divides
all four entries, while every common divisor of all four
entries necessarily divides both \(x\) and \(y\).

It follows from (3.4)--(3.6) that the Smith normal form of
\(M_\alpha\) is
\[
\begin{pmatrix}
g&0\\
0&n/g
\end{pmatrix}.
\tag{3.7}
\]

The cokernel of \(M_\alpha\) is precisely
\[
\OO_D/\alpha\OO_D.
\]
Since unimodular row and column operations preserve the
isomorphism class of the cokernel, (3.7) gives
\[
\OO_D/\alpha\OO_D
\simeq
\ZZ/g\ZZ\times\ZZ/(n/g)\ZZ.
\]
This proves (3.2).

Finally, the divisibility relation between the invariant
factors implies
\[
g\mid \frac ng,
\]
and hence \(g^2\mid n\). In particular,
\[
g\leq\sqrt n.
\]
The two invariant factors are equal precisely when
\[
g=\frac ng,
\]
which is equivalent to \(g^2=n\). Therefore they are unequal
if and only if \(g<\sqrt n\).
\end{proof}

\begin{remark}[Interpretation of the invariant factors]
The decomposition in (3.2) gives the following alternatives:
\[
\begin{array}{c|c}
\text{Condition on }g&\text{Structure of the group}\\ \hline
g=1&
\ZZ/n\ZZ\quad\text{(cyclic)}\\[2mm]
1<g<\sqrt n&
\ZZ/g\ZZ\times\ZZ/(n/g)\ZZ
\quad\text{(noncyclic, unequal factors)}\\[2mm]
g=\sqrt n&
\ZZ/\sqrt n\ZZ\times\ZZ/\sqrt n\ZZ
\quad\text{(equal factors).}
\end{array}
\]
\end{remark}

The preceding lemma concerns the quotient
\(\OO_D/\alpha\OO_D\). We next explain why the same
invariant factors occur for the kernel of the corresponding
endomorphism on every elliptic curve with endomorphism ring
\(\OO_D\).

\begin{proposition}[Independence of the ideal class]
\label{prop:ideal-quotient}
Let \(I\) be a proper invertible fractional
\(\OO_D\)-ideal and let
\[
E=\CC/I.
\]
For \(\alpha\in\OO_D\setminus\{0\}\), the finite abelian
groups
\[
I/\alpha I
\qquad\text{and}\qquad
\OO_D/\alpha\OO_D
\]
are isomorphic.

Consequently, if
\[
\alpha=x+y\omega_D,\qquad
\Norm(\alpha)=n,\qquad
g=\gcd(x,y),
\]
then the kernel of the endomorphism induced by multiplication
by \(\alpha\) on \(E\) has invariant factors
\[
g,\qquad n/g.
\]
\end{proposition}

\begin{proof}
Let \(\ell\) be a rational prime. Since \(I\) is an invertible
\(\OO_D\)-ideal, its localization at \(\ell\) is an invertible
module over the semilocal ring
\[
\OO_D\otimes_{\ZZ}\ZZ_\ell.
\]
An invertible module over a semilocal ring is free of rank one.
Thus
\[
I\otimes_{\ZZ}\ZZ_\ell
\simeq
\OO_D\otimes_{\ZZ}\ZZ_\ell
\]
as modules over \(\OO_D\otimes_{\ZZ}\ZZ_\ell\).

Under such an isomorphism, multiplication by \(\alpha\)
corresponds to multiplication by \(\alpha\). Hence
\[
(I/\alpha I)\otimes_{\ZZ}\ZZ_\ell
\simeq
(\OO_D/\alpha\OO_D)\otimes_{\ZZ}\ZZ_\ell.
\tag{3.8}
\]

Both quotients are finite abelian groups. A finite abelian
group is determined, up to isomorphism, by its primary
components. Since (3.8) holds for every rational prime
\(\ell\), the two finite abelian groups are isomorphic.

Finally, multiplication by \(\alpha\) on
\[
E=\CC/I
\]
is the endomorphism
\[
[z]\longmapsto[\alpha z].
\]
Its kernel is
\[
\ker(m_\alpha)
=
\alpha^{-1}I/I.
\]
Multiplication by \(\alpha\) induces an isomorphism
\[
\alpha^{-1}I/I
\longrightarrow
I/\alpha I.
\]
Therefore the kernel has the same invariant factors as
\(\OO_D/\alpha\OO_D\). The conclusion follows from
Lemma~\ref{lem:invariant-factors}.
\end{proof}

\begin{remark}
Proposition~\ref{prop:ideal-quotient} is the reason that the
invariant-factor calculation does not depend on the ideal
class representing the elliptic curve. Once an admissible
endomorphism exists for a given order \(\OO_D\), the same
kernel type occurs for every isomorphism class of elliptic
curves with endomorphism ring \(\OO_D\).
\end{remark}

\section{Admissible discriminants and the counting theorem}

We now formulate the arithmetic description of the counting
problem.

Throughout this section, \(n>1\) is a fixed positive integer.

\begin{definition}[Admissible discriminant sets]
\label{def:admissible}
Let \(D<0\) be a quadratic-order discriminant. We say that
\(D\) is cyclic-admissible for \(n\) if there exist
\(x,y\in\ZZ\), with \(y\ne0\), such that
\[
Q_D(x,y)=n,
\qquad
\gcd(x,y)=1.
\]
We define
\[
\DD_{\mathrm{cyc}}(n)
=
\left\{
\begin{array}{l|l}
D<0,\ D\equiv0,1\pmod4&
\begin{array}{l}
\exists\,x,y\in\ZZ,\ y\ne0,\\
Q_D(x,y)=n,\ \gcd(x,y)=1
\end{array}
\end{array}
\right\}.
\tag{4.1}
\]

Similarly, we define the noncyclic-admissible set by
\[
\DD_{\mathrm{nc}}(n)
=
\left\{
\begin{array}{l|l}
D<0,\ D\equiv0,1\pmod4&
\begin{array}{l}
\exists\,x,y\in\ZZ,\ y\ne0,\\
Q_D(x,y)=n,\\
1<\gcd(x,y)<\sqrt n
\end{array}
\end{array}
\right\}.
\tag{4.2}
\]

Finally, put
\[
\DD_{\mathrm{gen}}(n)
=
\DD_{\mathrm{cyc}}(n)
\cup
\DD_{\mathrm{nc}}(n).
\tag{4.3}
\]
Equivalently,
\[
\DD_{\mathrm{gen}}(n)
=
\left\{
\begin{array}{l|l}
D<0,\ D\equiv0,1\pmod4&
\begin{array}{l}
\exists\,x,y\in\ZZ,\ y\ne0,\\
Q_D(x,y)=n,\\
1\leq\gcd(x,y)<\sqrt n
\end{array}
\end{array}
\right\}.
\tag{4.4}
\]
\end{definition}

For clarity, the general admissible set in (4.3) concerns
kernels with unequal invariant factors. The equal-factor
case is excluded.

\begin{definition}[Numbers of elliptic-curve classes]
Let
\[
N_{\mathrm{cyc}}(n)
\]
denote the number of isomorphism classes of complex elliptic
curves \(E\) for which there exists a cyclic subgroup
\(C\subset E\) of order \(n\) such that \(E/C\simeq E\).

Let
\[
N_{\mathrm{nc}}(n)
\]
denote the number of isomorphism classes of complex elliptic
curves \(E\) for which there exists a noncyclic subgroup
\(C\subset E\) of order \(n\), with unequal invariant factors,
such that \(E/C\simeq E\).

Finally, let
\[
N_{\mathrm{gen}}(n)
\]
denote the number of isomorphism classes of complex elliptic
curves \(E\) for which there exists a subgroup \(C\subset E\)
of order \(n\), with unequal invariant factors, such that
\(E/C\simeq E\), without imposing cyclicity.
\end{definition}

We can now state the main result.

\begin{theorem}[Class-number reformulation and counting theorem]
\label{thm:main}
Let \(n>1\). With the notation of
Definition~\ref{def:admissible}, the numbers of isomorphism
classes of complex elliptic curves in the cyclic,
noncyclic, and general cases are given respectively by
\[
\boxed{
N_{\mathrm{cyc}}(n)
=
\sum_{D\in\DD_{\mathrm{cyc}}(n)}h(D).
}
\tag{4.5}
\]
Furthermore,
\[
\boxed{
N_{\mathrm{nc}}(n)
=
\sum_{D\in\DD_{\mathrm{nc}}(n)}h(D).
}
\tag{4.6}
\]
Finally,
\[
\boxed{
N_{\mathrm{gen}}(n)
=
\sum_{D\in\DD_{\mathrm{gen}}(n)}h(D).
}
\tag{4.7}
\]

In particular, the number of elliptic-curve classes in the
general case is the sum of the proper ideal class numbers of
the distinct admissible imaginary quadratic orders.
\end{theorem}

\begin{proof}
We first recall the relation between the quotient condition
and endomorphisms.

Suppose that \(E\) is a complex elliptic curve and that
\(C\subset E\) is a finite subgroup of order \(n\) satisfying
\[
E/C\simeq E.
\]
Let
\[
\pi:E\longrightarrow E/C
\]
be the natural quotient isogeny. Choose an isomorphism
\[
\psi:E/C\longrightarrow E.
\]
Then
\[
\alpha=\psi\circ\pi:E\longrightarrow E
\]
is an endomorphism of \(E\), and it satisfies
\[
\deg(\alpha)=n,
\qquad
\ker(\alpha)=C.
\tag{4.8}
\]

Conversely, if \(\alpha\in\End(E)\) is an endomorphism of
degree \(n\), then
\[
E/\ker(\alpha)\simeq E.
\tag{4.9}
\]
Thus the quotient condition is equivalent to the existence
of an endomorphism of degree \(n\) with the prescribed
kernel type.

We now consider the cases in the statement.

\medskip
\noindent
\textbf{Step 1: The endomorphism ring is an imaginary
quadratic order.}

If an endomorphism of a complex elliptic curve is an integer
multiplication, say \([m]\), then its kernel is
\[
E[m]\simeq
\ZZ/m\ZZ\times\ZZ/m\ZZ.
\]
Its invariant factors are equal.

Therefore, if the kernel has unequal invariant factors, the
endomorphism is not an integer multiplication. In this
situation \(E\) has complex multiplication, and its
endomorphism ring is an imaginary quadratic order
\(\OO_D\) for some discriminant \(D<0\).

Conversely, for an elliptic curve with endomorphism ring
\(\OO_D\), every endomorphism is represented by an element
of \(\OO_D\).

\medskip
\noindent
\textbf{Step 2: The norm equation.}

Write
\[
\OO_D=\ZZ[\omega_D],
\qquad
\omega_D=\frac{\varepsilon_D+\sqrt D}{2}.
\]
Every endomorphism \(\alpha\in\OO_D\) has the form
\[
\alpha=x+y\omega_D,
\qquad x,y\in\ZZ.
\]
Its degree is its norm:
\[
\deg(\alpha)=\Norm(\alpha)
=Q_D(x,y).
\]
Thus the degree condition \(\deg(\alpha)=n\) becomes
\[
Q_D(x,y)=n.
\tag{4.10}
\]

By Lemma~\ref{lem:invariant-factors} and
Proposition~\ref{prop:ideal-quotient}, the kernel of this
endomorphism has invariant factors
\[
g,\qquad n/g,
\qquad
g=\gcd(x,y).
\tag{4.11}
\]

In particular, the kernel is cyclic if and only if \(g=1\).
It is noncyclic with unequal invariant factors if and only if
\[
1<g<\sqrt n.
\]
The kernel has unequal invariant factors, without imposing
cyclicity, if and only if
\[
1\leq g<\sqrt n.
\]

These are exactly the conditions in
Definition~\ref{def:admissible}.

\medskip
\noindent
\textbf{Step 3: Counting the elliptic curves for a fixed
discriminant.}

Fix an admissible discriminant \(D\).

By Proposition~\ref{prop:CM}, the isomorphism classes of
complex elliptic curves with endomorphism ring
\(\OO_D\) are parametrized by the proper ideal classes of
\(\OO_D\). Their number is exactly \(h(D)\).

Moreover, if \(D\) is admissible, there exists an element
\[
\alpha=x+y\omega_D\in\OO_D
\]
of norm \(n\) satisfying the required greatest-common-divisor
condition.

For every proper ideal class represented by an invertible
ideal \(I\), multiplication by \(\alpha\) induces an
endomorphism of
\[
E_I=\CC/I.
\]
By Proposition~\ref{prop:ideal-quotient}, the kernel of this
endomorphism has the same invariant factors as
\[
\OO_D/\alpha\OO_D.
\]
Hence every isomorphism class of elliptic curves with
endomorphism ring \(\OO_D\) has the required kernel type.

Conversely, if an elliptic curve with endomorphism ring
\(\OO_D\) admits a subgroup of order \(n\) with the required
quotient and kernel type, the corresponding endomorphism
gives a solution of (4.10) satisfying the corresponding
greatest-common-divisor condition. Thus \(D\) belongs to the
appropriate admissible set.

Consequently, for each admissible discriminant \(D\), all
\(h(D)\) isomorphism classes associated with \(\OO_D\)
contribute to the corresponding count; if \(D\) is not
admissible, none of them contributes.

\medskip
\noindent
\textbf{Step 4: Summation over the discriminants.}

Different endomorphism rings have different discriminants.
An elliptic curve cannot have two distinct endomorphism
rings as subrings of \(\CC\), up to the canonical
identification of its endomorphism ring with an order.

Therefore, the classes associated with distinct
discriminants are disjoint.

Summing the contribution \(h(D)\) over the cyclic-admissible
discriminants gives
\[
N_{\mathrm{cyc}}(n)
=
\sum_{D\in\DD_{\mathrm{cyc}}(n)}h(D).
\]
The same argument gives
\[
N_{\mathrm{nc}}(n)
=
\sum_{D\in\DD_{\mathrm{nc}}(n)}h(D)
\]
and
\[
N_{\mathrm{gen}}(n)
=
\sum_{D\in\DD_{\mathrm{gen}}(n)}h(D).
\]
This proves (4.5)--(4.7).
\end{proof}


\subsection{The equal-invariant-factor case}

We now consider the case excluded from Theorem \ref{thm:main}, namely the
case in which the two invariant factors of the subgroup are equal.
This case has a simple description and does not require the
complex multiplication hypothesis.

\begin{proposition}
Let \(n>1\) be an integer. There exists a finite subgroup
\(C\subset E(\mathbb C)\) of order \(n\), with equal invariant
factors, satisfying
\[
E/C\simeq E,
\]
if and only if \(n=m^2\) for some integer \(m>1\).
In this case, for every complex elliptic curve \(E\), the
unique such subgroup is
\[
C=E[m].
\]
Consequently, if \(n\) is a perfect square, every complex
elliptic curve satisfies the required condition.
\end{proposition}

\begin{proof}
Suppose first that the invariant factors of \(C\) are equal.
Then, for some positive integer \(m\),
\[
C\simeq
\mathbb Z/m\mathbb Z\times
\mathbb Z/m\mathbb Z.
\]
In particular,
\[
n=|C|=m^2.
\]
Moreover, every element of \(C\) is annihilated by \(m\).
Thus
\[
C\subseteq E[m].
\]
On the other hand, for every complex elliptic curve,
\[
E[m]\simeq
\mathbb Z/m\mathbb Z\times
\mathbb Z/m\mathbb Z,
\]
and hence \(|E[m]|=m^2=|C|\). Therefore
\[
C=E[m].
\]

Conversely, suppose that \(n=m^2\), with \(m>1\).
For every complex elliptic curve \(E\), the multiplication
map
\[
[m]:E\longrightarrow E
\]
is an isogeny of degree \(m^2\), whose kernel is
\[
\ker([m])=E[m]\simeq
\mathbb Z/m\mathbb Z\times
\mathbb Z/m\mathbb Z.
\]
The first isomorphism theorem for groups applied to  to this surjective isogeny, the induced map $E/E[m] \to E$ is an isomorphism of complex elliptic curves, so
\[
E/E[m]\simeq E.
\]
Thus \(C=E[m]\) has all the required properties.

Finally, the inclusion \(C\subseteq E[m]\), together with
the equality of their orders, shows that no other subgroup
with equal invariant factors can occur.
\end{proof}

\begin{remark}
The equal-invariant-factor case is therefore completely
determined by multiplication maps. If \(n\) is not a perfect
square, no subgroup of order \(n\) can have equal invariant
factors. If \(n=m^2\), every complex elliptic curve admits
the subgroup \(E[m]\), and this subgroup is uniquely
determined.

In particular, when \(n\) is a perfect square, the number
of isomorphism classes of complex elliptic curves satisfying
the equal-invariant-factor condition is uncountably infinite.
Indeed, every complex elliptic curve satisfies it, and
the set of complex elliptic-curve isomorphism classes is
parametrized by the \(j\)-invariant in \(\mathbb C\).

This explains why the equal-invariant-factor case is not
included in the finite class-number sums of Theorem \ref{thm:main}:
it is governed by multiplication maps and includes
elliptic curves without complex multiplication.
\end{remark}

\begin{remark}
The same conclusion is consistent with the norm-form
description used in Section 3. Indeed, if
\[
\alpha=x+y\omega_D,\qquad
N(\alpha)=n,\qquad
g=\gcd(x,y),
\]
then equality of the invariant factors is equivalent to
\[
g=\frac{n}{g},
\qquad\text{or equivalently}\qquad n=g^2.
\]
Writing \(x=gu\) and \(y=gv\), with
\(\gcd(u,v)=1\), gives
\[
Q_D(u,v)=1.
\]
Thus \(\beta=u+v\omega_D\) is a unit of
\(\mathcal O_D\), and
\[
\alpha=g\beta.
\]
Consequently, the corresponding endomorphism is multiplication
by \(g\), followed by an automorphism of the elliptic curve.
Its kernel is precisely \(E[g]\).
\end{remark}

Note that the equal-invariant-factor case is not just another finite class-number contribution. It is a separate, completely solved case, and it contains all complex elliptic curves whenever \(n\) is a square.

\section{Finiteness and arithmetic bounds}

The class-number formulas in Theorem~\ref{thm:main} involve
only finitely many discriminants. We establish an explicit
bound.

\begin{proposition}[Bound for admissible discriminants]
\label{prop:bound}
Let \(n>1\), and suppose that \(D\) is admissible for \(n\).
Then
\[
\boxed{-4n\leq D\leq-3.}
\tag{5.1}
\]
In particular, every admissible discriminant belongs to the
finite set
\[
\left\{
D\in\ZZ:
-4n\leq D\leq-3,\quad
D\equiv0\text{ or }1\pmod4
\right\}.
\tag{5.2}
\]
\end{proposition}

\begin{proof}
By admissibility, there exist integers \(x,y\), with
\(y\ne0\), such that
\[
x^2+\varepsilon_Dxy+q_Dy^2=n.
\]
Multiplying by \(4\) and completing the square, we obtain
\[
\begin{aligned}
4n
&=4x^2+4\varepsilon_Dxy+4q_Dy^2\\
&=(2x+\varepsilon_Dy)^2
  +(4q_D-\varepsilon_D^2)y^2\\
&=(2x+\varepsilon_Dy)^2+|D|y^2.
\end{aligned}
\tag{5.3}
\]
Since \(y\ne0\), we have \(y^2\geq1\), and therefore
\[
4n\geq |D|.
\]
Thus
\[
D\geq-4n.
\]

On the other hand, \(D\) is a negative quadratic-order
discriminant, so
\[
D\equiv0\ \text{or}\ 1\pmod4.
\]
The largest negative integer satisfying one of these
congruences is \(-3\). Hence
\[
D\leq-3.
\]
Combining the two inequalities proves (5.1).
\end{proof}

\begin{remark}[A divisibility condition]
If
\[
g=\gcd(x,y),
\]
write
\[
x=gu,\qquad y=gv,
\qquad \gcd(u,v)=1.
\]
Since \(Q_D\) is homogeneous of degree two,
\[
n=Q_D(x,y)=g^2Q_D(u,v).
\]
Thus
\[
\boxed{g^2\mid n.}
\tag{5.4}
\]
This gives an additional restriction when enumerating
noncyclic admissible discriminants.
\end{remark}

\begin{remark}[A finite enumeration procedure]
Proposition~\ref{prop:bound} gives a direct finite procedure
for determining the admissible discriminants.

For a fixed \(n\), it is enough to examine the integers
\(D\) satisfying
\[
-4n\leq D\leq-3,
\qquad
D\equiv0\ \text{or}\ 1\pmod4.
\]
For each such \(D\), one tests the integral solutions of
\[
x^2+\varepsilon_Dxy+q_Dy^2=n,
\qquad y\ne0,
\]
and records the possible values of \(\gcd(x,y)\).

The discriminant \(D\) is included in the appropriate
admissible set if at least one solution satisfies the
corresponding condition. It is included only once, even
when several pairs \((x,y)\) satisfy the norm equation.
\end{remark}

\section{The overlap between the cyclic and noncyclic cases}

It is important to distinguish the admissible discriminant
sets from the corresponding sets of elliptic curves.

Although a particular subgroup is either cyclic or
noncyclic, a single elliptic curve may admit two different
subgroups of order \(n\), one cyclic and the other
noncyclic, each satisfying the quotient condition.

Consequently, the discriminant sets
\[
\DD_{\mathrm{cyc}}(n)
\quad\text{and}\quad
\DD_{\mathrm{nc}}(n)
\]
need not be disjoint.

By definition,
\[
\DD_{\mathrm{gen}}(n)
=
\DD_{\mathrm{cyc}}(n)\cup
\DD_{\mathrm{nc}}(n).
\]
Therefore, inclusion-exclusion gives the identity
\[
\begin{aligned}
N_{\mathrm{gen}}(n)
={}&N_{\mathrm{cyc}}(n)+N_{\mathrm{nc}}(n)\\
&-
\sum_{D\in
\DD_{\mathrm{cyc}}(n)\cap
\DD_{\mathrm{nc}}(n)}
h(D).
\end{aligned}
\tag{6.1}
\]

In particular, the general count is not necessarily the sum
of the cyclic and noncyclic counts.

\begin{remark}
The intersection in (6.1) has a precise interpretation:
for every discriminant in the intersection, every elliptic
curve with the corresponding endomorphism ring admits both
a cyclic kernel and a noncyclic kernel of order \(n\) with
the required quotient property.
\end{remark}

\section{Numerical examples}

We illustrate Theorem~\ref{thm:main} with four values of
\(n\). The computations use the finite discriminant
enumeration described in Section 5, together with the
class-number sums in (4.5)--(4.7).

The resulting counts are displayed in Table~\ref{tab:counts}.

\begin{table}[ht]
\centering
\caption{Class-number counts for selected values of \(n\).}
\label{tab:counts}
\begin{tabular}{c|r|r|r|r}
\hline
\(n\)
&
\(N_{\mathrm{cyc}}(n)\)
&
\(N_{\mathrm{nc}}(n)\)
&
\(N_{\mathrm{gen}}(n)\)
&
Overlap
\\
\hline
150  & 340  & 10  & 348  & 2\\
297  & 409  & 43  & 440  & 12\\
1571 & 1605 & 0   & 1605 & 0\\
2012 & 3039 & 524 & 3563 & 0\\
\hline
\end{tabular}
\end{table}

Here the overlap column denotes
\[
\sum_{D\in
\DD_{\mathrm{cyc}}(n)\cap
\DD_{\mathrm{nc}}(n)}
h(D).
\]
Thus, for example, when \(n=150\), the general count is
\[
N_{\mathrm{gen}}(150)
=
340+10-2
=
348.
\]
Similarly,
\[
N_{\mathrm{gen}}(297)
=
409+43-12
=
440.
\]

For \(n=1571\), there are no noncyclic admissible
discriminants, so the general count coincides with the
cyclic count. For \(n=2012\), the cyclic and noncyclic
contributions are disjoint, and hence
\[
N_{\mathrm{gen}}(2012)
=
3039+524
=
3563.
\]

These examples illustrate why the overlap correction in
(6.1) is necessary. The Magma code used for these computations is available from the author upon request.

\section{Concluding remarks}

The main purpose of this paper has been to express the
counting problem for complex elliptic curves admitting
self-quotients in terms of class numbers of imaginary
quadratic orders.

The underlying classification of elliptic curves admitting
subgroups $C$ with
\[
E/C\simeq E
\]
is part of the established theory of complex multiplication
and was treated algorithmically in earlier work [1, 2, 3].
The contribution presented here is an arithmetic
reformulation of the counting problem, together with
explicit class-number formulas and a finite procedure
for their computation.

For an imaginary quadratic order $\mathcal O_D$, the
relevant endomorphisms are described by the norm equation
\[
x^2+\varepsilon_Dxy+q_Dy^2=n.
\]
The Smith normal form of the multiplication matrix shows
that the invariant factors of the kernel are
\[
\gcd(x,y),\qquad
\frac{n}{\gcd(x,y)}.
\]
This gives explicit arithmetic criteria for the cyclic,
noncyclic, and general cases with unequal invariant factors.

The CM correspondence then implies that every admissible
discriminant $D$ contributes exactly $h(D)$ isomorphism
classes of elliptic curves. We therefore obtain the
class-number identities
\[
N_{\mathrm{cyc}}(n)
=
\sum_{D\in\mathcal D_{\mathrm{cyc}}(n)}h(D),
\]
\[
N_{\mathrm{nc}}(n)
=
\sum_{D\in\mathcal D_{\mathrm{nc}}(n)}h(D),
\]
and
\[
N_{\mathrm{gen}}(n)
=
\sum_{D\in\mathcal D_{\mathrm{gen}}(n)}h(D).
\]
The discriminant bound
\[
-4n\leq D\leq -3
\]
makes these formulas finite and provides a direct
arithmetic procedure for their computation.

The relation between the cyclic and noncyclic cases also
requires care. Since a single elliptic curve may admit
both types of subgroups of the prescribed order, the
corresponding discriminant sets need not be disjoint.
The general count is therefore obtained by subtracting
the class-number contribution of their intersection.
The numerical examples illustrate the necessity of this
overlap correction.

Finally, we have completely determined the case of equal
invariant factors. Proposition 4.4 shows that such a
subgroup exists precisely when $n=m^2$, and that it is
uniquely given by $E[m]$ for every complex elliptic curve
$E$. Thus, when $n$ is a perfect square, every complex
elliptic curve satisfies the required condition, and the
number of isomorphism classes is uncountably infinite.
This case is therefore qualitatively different from the
unequal-invariant-factor cases: it is governed by
multiplication maps and includes elliptic curves without
complex multiplication.



\noindent
\textbf{Radu Gaba} Institute of
Mathematics ``Simion Stoilow'' of the Romanian Academy,
Bucharest and the National University of Science Politehnica Bucharest, Romania.
E-mail: \texttt{radu.gaba@imar.ro}.

\medskip
%
%
%

\end{document}